\pdfoutput=1
\RequirePackage{ifpdf}
\ifpdf 
\documentclass[pdftex]{sigma}
\else
\documentclass{sigma}
\fi

\usepackage{amscd,array,stmaryrd,mathrsfs,bbm}
\usepackage{colortbl} 
\definecolor{mygray}{gray}{.7}

\numberwithin{equation}{section}
\newtheorem{thm}{Theorem}[section]
\newtheorem*{MT}{Main Theorem}
\theoremstyle{definition}
\newtheorem{defi}[thm]{Definition}
\newtheorem{rem}[thm]{Remark}

\def\k{\mathbbm{k}}

\begin{document}

\allowdisplaybreaks

\newcommand{\arXivNumber}{1705.04913}

\renewcommand{\PaperNumber}{064}

\FirstPageHeading

\ShortArticleName{A Generalization of the Doubling Construction for Sums of Squares Identities}

\ArticleName{A Generalization of the Doubling Construction\\ for Sums of Squares Identities}

\Author{Chi ZHANG~$^\dag$ and Hua-Lin HUANG~$^\ddag$}

\AuthorNameForHeading{C.~Zhang and H.-L.~Huang}

\Address{$^\dag$~School of Mathematics, Shandong University, Jinan 250100, China}
\EmailD{\href{mailto:chizhang@mail.sdu.edu.cn}{chizhang@mail.sdu.edu.cn}}

\Address{$^\ddag$~School of Mathematical Sciences, Huaqiao University, Quanzhou 362021, China}
\EmailD{\href{mailto:hualin.huang@hqu.edu.cn}{hualin.huang@hqu.edu.cn}}

\ArticleDates{Received May 16, 2017, in f\/inal form August 13, 2017; Published online August 16, 2017}

\Abstract{The doubling construction is a fast and important way to generate new solutions to the Hurwitz problem on sums of squares identities from any known ones. In this short note, we generalize the doubling construction and obtain from any given admissible triple $[r,s,n]$ a series of new ones $[r+\rho(2^{m-1}),2^ms,2^mn]$ for all positive integer $m$, where $\rho$ is the Hurwitz--Radon function.}


\Classification{11E25}

\section{Introduction}

In his seminal paper \cite{H}, Hurwitz addressed the famous problem: Determine all the sums of squares identities
\begin{gather}\label{sqi}
\big(x_1^2+x_2^2+\cdots +x_r^2\big)\big(y_1^2+y_2^2+\cdots+y_s^2\big)= z_1^2+z_2^2+\cdots+z_n^2,
\end{gather}
where $X=(x_1,x_2,\dots,x_r)$ and $Y=(y_1,y_2,\dots,y_s)$ are systems of indeterminants and every $z_k$ is a bilinear form of $X$ and $Y$ with coef\/f\/icients in some given f\/ield. If there does exist such an identity, we call $[r,s,n]$ an \emph{admissible triple}. This problem of Hurwitz has close connections to various topics in algebra, arithmetic, combinatorics, geometry, topology, etc. Many mathematicians have studied this Hurwitz problem during the past century. See~\cite{S} for an overview.

Though a complete solution to the Hurwitz problem is still far out of reach at present, many admissible triples have been obtained in the literature. In particular, the admissible triples of form $[r, n, n]$ was settled independently by Hurwitz in~\cite{H1} and by Radon in~\cite{Rad}. The celebrated Hurwitz--Radon theorem states that $[r, n, n]$ is admissible if and only if $r \le \rho(n)$ where $\rho$ is the Hurwitz--Radon function def\/ined by $\rho(n)=8\alpha+2^\beta$ if $n=2^{4\alpha+\beta}(2\gamma+1)$ with $0 \le \beta \le 3$. In the early 1980s, Yuzvinsky introduced the novel idea of orthogonal pairings~\cite{Y1} and proposed in~\cite{Y2} the following three families of admissible triples in the neighborhood of the Hurwitz--Radon triples
\begin{gather}\label{eq:Yuzvinsky}
\big[2n+2, 2^n-\varphi(n), 2^n\big], \qquad \hbox{where} \quad
\varphi(n)= \begin{cases}
 \binom{n }{n/2}, & n \equiv 0 \mod 4, \\
 2\binom{n-1}{(n-1)/2}, & n \equiv 1 \mod 4, \\
 4\binom{n-2}{(n-2)/2}, & n \equiv 2 \mod 4.
 \end{cases}
\end{gather}
The f\/irst two families are conf\/irmed in \cite{LS} and the third in \cite{HHZ}. Moreover, some new families of admissible triples are constructed in \cite{HHZ, LMGO}.

Another natural idea of constructing admissible triples is to f\/ind some general procedures to generate new ones from known ones. Among which, the doubling construction, i.e., generating an admissible triple $[r+1,2s,2n]$ from any given triple $[r,s,n]$, is a fast and important program. In the present note we consider arbitrarily iterated doubling constructions and the aim is to optimize the obvious triples $[r+m,2^ms,2^mn]$. The well-known Hurwitz--Radon triples suggest that the f\/irst item might be properly increased as the form given below via the function~$\rho$.

\begin{MT} If $[r,s,n]$ is admissible, then so is $[r+\rho(2^{m-1}),2^ms,2^mn]$ for all positive integer $m$.
\end{MT}

Though our observation arises from the idea used in \cite{HHZ}, it turns out that a more elementary approach of matrices will suf\/f\/ice for a proof.

\section{Proof of the main theorem}

First, we introduce the so-called admissible matrices to reformulate the Hurwitz problem. Then, as a trial we provide a proof via admissible matrices for the classical doubling construction. Finally, we extend the idea to iterated doubling constructions and complete the proof for the main theorem.

\subsection{Admissible matrices}

The notion of admissible matrices arises naturally from an attempt to reformulate the sums of squares identity (\ref{sqi}) by a system of polynomial equations. Indeed, in (\ref{sqi}) if we write
\begin{gather*}
z_k = \sum_{\substack{1 \leq i \leq r \\ 1 \leq j \leq s}} c_{i,j,k} x_iy_j,
\end{gather*}
then it is easy to see that the identity (\ref{sqi}) is equivalent to the following system of algebraic equations
\begin{alignat}{3}
 & \sum\limits_{k=1}^n c_{i,j,k}^2 = 1, \qquad & & 1 \leq i \leq r , \quad 1 \leq j \leq s, &\nonumber\\
 & \sum\limits_{k=1}^n c_{i_1,j,k}c_{i_2,j,k} = 0,\qquad  & & 1 \leq i_1 < i_2 \leq r , \quad 1 \le j \le s, & \nonumber\\
 & \sum\limits_{k=1}^n c_{i,j_1,k}c_{i,j_2,k} = 0,\qquad  & & 1 \leq i \leq r , \quad 1 \leq j_1 < j_2 \leq s, &\nonumber\\
 & \sum\limits_{k=1}^n (c_{i_1,j_1,k}c_{i_2,j_2,k}+c_{i_1,j_2,k}c_{i_2,j_1,k}) = 0,\qquad  & & 1 \leq i_1 < i_2 \leq r, \quad 1 \leq j_1 < j_2 \leq s.& \label{primary conditions}
 \end{alignat}
In the rest of the note, we always regard the resulting cuboid $A:=(c_{i,j,k})_{r \times s \times n}$ as an $r \times s$ matrix with $(i,j)$-entry the $n$-dimensional vector $A_{i,j}:=(c_{i,j,1}, c_{i,j,2}, \dots, c_{i,j,n})$. Taking the formal inner product on $n$-dimensional vectors, namely $\langle (u_1, u_2, \dots, u_n), (v_1, v_2, \dots, v_n) \rangle:=u_1v_1 + u_2v_2 + \cdots +u_nv_n$, then~(\ref{primary conditions}) can be rewritten as
\begin{alignat}{4}
& (1) \ \ && \langle A_{i,j}, A_{i,j}\rangle = 1, \qquad  & & 1 \leq i \leq r , \quad 1 \leq j \leq s, &\nonumber\\
& (2) \ \ & & \langle A_{i_1,j}, A_{i_2,j}\rangle = 0, \qquad  & & 1 \leq i_1 < i_2 \leq r , \quad 1 \leq j\leq s, &\nonumber\\
& (3) \ \ && \langle A_{i,j_1}, A_{i,j_2}\rangle = 0,\qquad & & 1 \leq i \leq r , \quad 1 \leq j_1 < j_2 \leq s, &\nonumber\\
& (4) \ \ && \langle A_{i_1,j_1}, A_{i_2,j_2}\rangle + \langle A_{i_1,j_2}, A_{i_2,j_1}\rangle = 0,\qquad & & 1 \leq i_1 < i_2 \leq r, \quad 1 \leq j_1 < j_2 \leq s.&\label{eq:admissible matrix}
 \end{alignat}
Obviously, the existence of such a matrix $A$ is equivalent to the existence of an admissible triple of size $[r,s,n]$. In keeping the terminologies coherent, such matrices are said to be \emph{admissible}.

\subsection{The doubling construction revisited}

For a better explanation of our method, f\/irstly we provide a proof by admissible matrices for the classical doubling construction. Some preparing def\/initions and notations are necessary. Let $\k$ be a f\/ield of characteristic not~2.

\begin{defi}
Fix two integers $n$ and $m$. A vector in $\alpha=(\alpha_1, \alpha_2, \dots, \alpha_{mn}) \in \k^{mn}$ is said to be in {\it level} $k \in \{ 1, 2, \dots, m\} $ if its arguments $\alpha_l$ are $0$ unless $(k-1)n+1 \le l \le kn$. Let $\beta \in \k^{n}$. We call $\gamma \in \k^{mn}$ a~{\it positive} copy of $\beta$ in level $k$ if $\gamma_{(k-1)n+i}=\beta_{i}$ $(1 \leq i \leq n)$ and other arguments of~$\gamma$ are~$0$. Similarly, we call $\gamma$ a~{\it negative} copy of $\beta$ in level $k$ if $\gamma_{(k-1)n+i}=-\beta_{i}$ $(1 \leq i \leq n)$ and other arguments of $\gamma$ are~$0$.
\end{defi}
\begin{rem}\label{rem:Euclid product and levels} Let $\alpha_1, \alpha_2 \in \k^{mn}$ and $\beta_1, \beta_2 \in \k^n$.
\begin{enumerate}\itemsep=0pt
 \item If $\alpha_1$ is a copy of $\beta_1$ in level $k$, $\alpha_2$ is a copy of $\beta_2$ in level $l$ and they have the same sign, then $\langle \alpha_1, \alpha_2 \rangle = \langle \beta_1, \beta_2 \rangle$ if $k=l$, $\langle \alpha_1, \alpha_2 \rangle = 0$ if $k \neq l$.
 \item If $\alpha_1$ is a copy of $\beta_1$ in level $k$, $\alpha_2$ is a copy of $\beta_2$ in level $l$ and they have dif\/ferent signs, then $\langle \alpha_1, \alpha_2 \rangle = -\langle \beta_1, \beta_2 \rangle$ if $k=l$, $\langle \alpha_1, \alpha_2 \rangle = 0$ if $k \neq l$.
\end{enumerate}
\end{rem}

Given an admissible triple of size $[r,s,n]$, we have a corresponding $r \times s$ admissible mat\-rix~$A$. We shall construct an $(r+1) \times 2s$ admissible matrix $B$ whose entries are $2n$-dimensional vectors as follows:
\begin{enumerate}\itemsep=0pt
\item[1)] $B_{i,j}$ is a positive copy of $A_{i,j}$ for $1 \leq i \leq r$, $1 \leq j \leq s$ in level $1$,
\item[2)] $B_{i,s+j}$ is a positive copy of $A_{i,j}$ for $2 \leq i \leq r$, $1 \leq j \leq s$ in level $2$,
\item[3)] $B_{r+1,j}$ is a positive copy of $A_{1,j}$ for $1 \leq j \leq s$ in level $2$,
\item[4)] $B_{r+1,s+j}$ is a positive copy of $A_{1,j}$ for $1 \leq j \leq s$ in level $1$,
\item[5)] $B_{1,s+j}$ is a negative copy of $A_{1,j}$ for $1 \leq j \leq s$ in level $2$.
\end{enumerate}

We give a detailed verif\/ication of the admissibility of $B$ and hope this will shed some light on the study of iterated doubling constructions.
\begin{enumerate}\itemsep=0pt
\item Every $B_{i,j}$ is a copy of some entry $A_{k,l}$ of $A$ , so $\langle B_{i,j}, B_{i,j} \rangle = \langle A_{k,l}, A_{k,l} \rangle = 1$, hence~(1) of~(\ref{eq:admissible matrix}) holds.
\item $\langle B_{i,j_1}, B_{i,j_2} \rangle = 0$ $(j_1 < j_2)$. Indeed, if $1 \leq j_1 \leq s$ and $s+1 \leq j_2 \leq 2s$, then the two vectors are in dif\/ferent levels; if $1 \leq i \leq r+1$, $1 \leq j_1 < j_2 \leq s$, $\langle B_{i,j_1}, B_{i,j_2} \rangle = \langle A_{i,j_1}, A_{i,j_2} \rangle = 0$. A~similar argument works for $s+1 \leq j_1 < j_2 \leq 2s$. So~(2) of~(\ref{eq:admissible matrix}) holds. In the same way, one can show that~(3) of~(\ref{eq:admissible matrix}) holds.
\item For (4) of (\ref{eq:admissible matrix}), we need to verify $\langle B_{i_1,j_1}, B_{i_2,j_2} \rangle + \langle B_{i_1,j_2}, B_{i_2,j_1} \rangle = 0$ $(i_1 < i_2$, $j_1 < j_2)$.
\begin{enumerate}\itemsep=0pt
\item If $1 \leq j_1 \leq s < j_2 \leq 2s$ and $1 \leq i_1 < i_2 \leq r$, $B_{i_1,j_1}$ and $B_{i_2,j_1}$ are in level $1$ and $B_{i_1,j_2}$ and $B_{i_2,j_2}$ are in level $2$. Hence the equation is obvious.
\item If $1 \leq j_1 \leq s$, $s+1 \leq j_2 \leq 2s$ and $1 \leq i_1 \leq r$, $i_2 = r+1$, $B_{i_1,j_1}$ and $B_{i_2,j_2}$ are in level $1$ and $B_{i_1,j_2}$ and $B_{i_2,j_1}$ are in level~$2$. $\langle B_{i_1,j_1}, B_{i_2,j_2} \rangle + \langle B_{i_1,j_2}, B_{i_2,j_1} \rangle = \langle A_{i_1,j_1}, A_{1,j_2} \rangle + \langle A_{i_1,j_2}, A_{1,j_1} \rangle = 0$. If $i_1 = 1$ and $j_2 = j_1 + s$, then $\langle B_{i_1,j_1}, B_{i_2,j_2} \rangle + \langle B_{i_1,j_2}, B_{i_2,j_1} \rangle = \langle A_{1,j_1}, A_{1,j_1} \rangle - \langle A_{1,j_1}, A_{1,j_1} \rangle = 1 - 1 = 0$.
\item If $1 \leq i_1 < i_2 \leq r$, $1 \leq j_1 < j_2 \leq s$ or $s+1 \leq j_1 < j_2 \leq 2s$, the four vectors are in the same level. If $i_1 = 1$, $s+1 \leq j_1\leq j_2 \leq 2s$, then $\langle B_{i_1,j_1}, B_{i_2,j_2} \rangle + \langle B_{i_1,j_2}, B_{i_2,j_1} \rangle = - \langle A_{i_1,j_1}, A_{i_2,j_2} \rangle - \langle A_{i_1,j_2}, A_{i_2,j_1} \rangle = 0$. Otherwise, $\langle B_{i_1,j_1}, B_{i_2,j_2} \rangle + \langle B_{i_1,j_2}, B_{i_2,j_1} \rangle = \langle A_{i_1,j_1}, A_{i_2,j_2} \rangle + \langle A_{i_1,j_2}, A_{i_2,j_1} \rangle = 0$.
\item If $1 \leq i_1 \leq r$, $i_2 = r+1$, $1 \leq j_1 < j_2 \leq s$ or $s+1 \leq j_1 < j_2 \leq 2s$, then~$B_{i_1,j_1}$ and~$B_{i_2,j_2}$ are in dif\/ferent levels and this is also the case for~$B_{i_1,j_2}$ and~$B_{i_2,j_1}$. Now the equation is clear.
\end{enumerate}
Thus, (4) of (\ref{eq:admissible matrix}) holds.
\end{enumerate}

\begin{rem}\label{discussion of judgement}
The matrix $B$ used in the above proof can be illustrated by the following table:
$$
\begin{tabular}{|m{100pt}<{\centering}|m{100pt}<{\centering}|}
\hline
1 & \cellcolor{mygray}{2} \\
\hline
 & \\
1 & 2\\
 & \\
\hline
2 & 1\\
\hline
\end{tabular}
$$
Here, cells in the f\/irst and the third rows are copies of the f\/irst row of $A$, and cells in the second rows are copies of the submatrix of $A$ obtained by deleting the f\/irst row. The number given in the center of a cell represents the level of the vectors therein. The sign of a cell is indicated by its color: white means positive, gray means negative.

Above all, the table provides a visual admissibility of $B$. The conditions (1)--(3) of~(\ref{eq:admissible matrix}) are immediate, as the vectors are either in dif\/ferent levels, or essentially can be considered within $A$. The same reasoning also works for~(4) of~(\ref{eq:admissible matrix}) in most cases. As for the case $i_1 = 1$, $i_2 = r+1$, $j_1 + s = j_2$, one further needs to take the signs into consideration. In fact, this also tells us in the very beginning how to manipulate the signs of the copies of cells so that~(4) holds. Of course, the signing is far from unique. Just for such $B$, we have 16 kinds of correct schemes as follows. These tables are useful in the following for the verif\/ication of the admissibility of bigger matrices.
$$
\begin{array}{cccc}
\begin{tabular}{|m{20pt}<{\centering}|m{20pt}<{\centering}|}
\hline
1 & \cellcolor{mygray}{2} \\
\hline
 & \\
1 & 2\\
 & \\
\hline
2 & 1\\
\hline
\end{tabular}
&
\begin{tabular}{|m{20pt}<{\centering}|m{20pt}<{\centering}|}
\hline
1 & \cellcolor{mygray}{2} \\
\hline
\cellcolor{mygray}{} & \cellcolor{mygray}{} \\
\cellcolor{mygray}{1} & \cellcolor{mygray}{2}\\
\cellcolor{mygray}{} & \cellcolor{mygray}{} \\
\hline
2 & 1\\
\hline
\end{tabular}
&
\begin{tabular}{|m{20pt}<{\centering}|m{20pt}<{\centering}|}
\hline
1 & \cellcolor{mygray}{2} \\
\hline
 & \\
1 & 2\\
 & \\
\hline
\cellcolor{mygray}{2} & \cellcolor{mygray}{1}\\
\hline
\end{tabular}
&
\begin{tabular}{|m{20pt}<{\centering}|m{20pt}<{\centering}|}
\hline
1 & \cellcolor{mygray}{2} \\
\hline
\cellcolor{mygray}{} & \cellcolor{mygray}{} \\
\cellcolor{mygray}{1} & \cellcolor{mygray}{2}\\
\cellcolor{mygray}{} & \cellcolor{mygray}{} \\
\hline
\cellcolor{mygray}{2} & \cellcolor{mygray}{1}\\
\hline
\end{tabular}
\end{array}
$$
$$
\begin{array}{cccc}
\begin{tabular}{|m{20pt}<{\centering}|m{20pt}<{\centering}|}
\hline
\cellcolor{mygray}{1} & 2 \\
\hline
 & \\
1 & 2\\
 & \\
\hline
2 & 1\\
\hline
\end{tabular}
&
\begin{tabular}{|m{20pt}<{\centering}|m{20pt}<{\centering}|}
\hline
\cellcolor{mygray}{1} & 2 \\
\hline
\cellcolor{mygray}{} & \cellcolor{mygray}{} \\
\cellcolor{mygray}{1} & \cellcolor{mygray}{2}\\
\cellcolor{mygray}{} & \cellcolor{mygray}{} \\
\hline
2 & 1\\
\hline
\end{tabular}
&
\begin{tabular}{|m{20pt}<{\centering}|m{20pt}<{\centering}|}
\hline
\cellcolor{mygray}{1} & 2 \\
\hline
 & \\
1 & 2\\
 & \\
\hline
\cellcolor{mygray}{2} & \cellcolor{mygray}{1}\\
\hline
\end{tabular}
&
\begin{tabular}{|m{20pt}<{\centering}|m{20pt}<{\centering}|}
\hline
\cellcolor{mygray}{1} & 2 \\
\hline
\cellcolor{mygray}{} & \cellcolor{mygray}{} \\
\cellcolor{mygray}{1} & \cellcolor{mygray}{2}\\
\cellcolor{mygray}{} & \cellcolor{mygray}{} \\
\hline
\cellcolor{mygray}{2} & \cellcolor{mygray}{1}\\
\hline
\end{tabular}
\end{array}
$$
$$
\begin{array}{cccc}
\begin{tabular}{|m{20pt}<{\centering}|m{20pt}<{\centering}|}
\hline
1 & 2 \\
\hline
\cellcolor{mygray}{} & \\
\cellcolor{mygray}{1} & 2\\
\cellcolor{mygray}{} & \\
\hline
\cellcolor{mygray}{2} & 1\\
\hline
\end{tabular}
&
\begin{tabular}{|m{20pt}<{\centering}|m{20pt}<{\centering}|}
\hline
1 & 2 \\
\hline
 & \cellcolor{mygray}{} \\
1 & \cellcolor{mygray}{2}\\
 & \cellcolor{mygray}{} \\
\hline
2 & \cellcolor{mygray}{1}\\
\hline
\end{tabular}
&
\begin{tabular}{|m{20pt}<{\centering}|m{20pt}<{\centering}|}
\hline
1 & 2 \\
\hline
\cellcolor{mygray}{} & \\
\cellcolor{mygray}{1} & 2\\
\cellcolor{mygray}{} & \\
\hline
2 & \cellcolor{mygray}{1}\\
\hline
\end{tabular}
&
\begin{tabular}{|m{20pt}<{\centering}|m{20pt}<{\centering}|}
\hline
1 & 2 \\
\hline
 & \cellcolor{mygray}{} \\
1 & \cellcolor{mygray}{2}\\
 & \cellcolor{mygray}{} \\
\hline
\cellcolor{mygray}{2} & 1\\
\hline
\end{tabular}
\end{array}
$$
$$
\begin{array}{cccc}
\begin{tabular}{|m{20pt}<{\centering}|m{20pt}<{\centering}|}
\hline
\cellcolor{mygray}{1} & \cellcolor{mygray}{2} \\
\hline
\cellcolor{mygray}{} & \\
\cellcolor{mygray}{1} & 2\\
\cellcolor{mygray}{} & \\
\hline
\cellcolor{mygray}{2} & 1\\
\hline
\end{tabular}
&
\begin{tabular}{|m{20pt}<{\centering}|m{20pt}<{\centering}|}
\hline
\cellcolor{mygray}{1} & \cellcolor{mygray}{2} \\
\hline
 & \cellcolor{mygray}{} \\
1 & \cellcolor{mygray}{2}\\
 & \cellcolor{mygray}{} \\
\hline
2 & \cellcolor{mygray}{1}\\
\hline
\end{tabular}
&
\begin{tabular}{|m{20pt}<{\centering}|m{20pt}<{\centering}|}
\hline
\cellcolor{mygray}{1} & \cellcolor{mygray}{2} \\
\hline
\cellcolor{mygray}{} & \\
\cellcolor{mygray}{1} & 2\\
\cellcolor{mygray}{} & \\
\hline
2 & \cellcolor{mygray}{1}\\
\hline
\end{tabular}
&
\begin{tabular}{|m{20pt}<{\centering}|m{20pt}<{\centering}|}
\hline
\cellcolor{mygray}{1} & \cellcolor{mygray}{2} \\
\hline
 & \cellcolor{mygray}{} \\
1 & \cellcolor{mygray}{2}\\
 & \cellcolor{mygray}{} \\
\hline
\cellcolor{mygray}{2} & 1\\
\hline
\end{tabular}
\end{array}
$$
\end{rem}

\subsection{Iterated doubling constructions}
Now we are ready to prove the main theorem. As before, let~$A$ be an admissible matrix corresponding to an admissible triple of size $[r,s,n]$. We will provide admissible matrices in terms of tables as in Remark~\ref{discussion of judgement} which will induce admissible triples of sizes $[r+2,4s,4n]$, $[r+4,8s,8n]$ and $[r+8,16s,16n]$.
\begin{center}
\begin{tabular}{|m{80pt}<{\centering}|m{80pt}<{\centering}|m{80pt}<{\centering}|m{80pt}<{\centering}|}
\hline
1 & \cellcolor{mygray}{2} & \cellcolor{mygray}{3} & 4\\
\hline
 & & & \\
1 & 2 & 3 & 4\\
 & & & \\
\hline
2 & 1 & \cellcolor{mygray}{4} & \cellcolor{mygray}{3} \\
\hline
3 & 4 & 1 & 2 \\
\hline
\end{tabular}\\
\centerline{the table of $[r+2,4s,4n]$}
\end{center}
\begin{center}
\begin{tabular}{|m{40pt}<{\centering}|m{40pt}<{\centering}|m{40pt}<{\centering}|m{40pt}<{\centering}|m{40pt}<{\centering}|m{40pt}<{\centering}|m{40pt}<{\centering}|m{40pt}<{\centering}|}
\hline
1 & \cellcolor{mygray}{2} & \cellcolor{mygray}{3} & 4 & \cellcolor{mygray}{5} & 6 & 7 & \cellcolor{mygray}{8}\\
\hline
 & & & & & & & \\
1 & 2 & 3 & 4 & 5 & 6 & 7 & 8\\
 & & & & & & & \\
\hline
2 & 1 & \cellcolor{mygray}{4} & \cellcolor{mygray}{3} & \cellcolor{mygray}{6} & \cellcolor{mygray}{5} & 8 & 7\\
\hline
3 & 4 & 1 & 2 & \cellcolor{mygray}{7} & \cellcolor{mygray}{8} & \cellcolor{mygray}{5} & \cellcolor{mygray}{6}\\
\hline
5 & 6 & 7 & 8 & 1 & 2 & 3 & 4\\
\hline
8 & \cellcolor{mygray}{7} & 6 & \cellcolor{mygray}{5} & \cellcolor{mygray}{4} & 3 & \cellcolor{mygray}{2} & 1\\
\hline
\end{tabular}\\
\centerline{the table of $[r+4,8s,8n]$}
\end{center}
\begin{center}
\begin{tabular}{|m{16pt}<{\centering}|m{16pt}<{\centering}|m{16pt}<{\centering}|m{16pt}<{\centering}|m{16pt}<{\centering}|m{16pt}<{\centering}|m{16pt}<{\centering}|m{16pt}<{\centering}|m{16pt}<{\centering}|m{16pt}<{\centering}|m{16pt}<{\centering}|m{16pt}<{\centering}|m{16pt}<{\centering}|m{16pt}<{\centering}|m{16pt}<{\centering}|m{16pt}<{\centering}|}
\hline
1 & \cellcolor{mygray}{2} & \cellcolor{mygray}{3} & 4 & \cellcolor{mygray}{5} & 6 & 7 & \cellcolor{mygray}{8} &
\cellcolor{mygray}{9} & 10 & 11 & \cellcolor{mygray}{12} & 13 & \cellcolor{mygray}{14} & \cellcolor{mygray}{15} & 16\\
\hline
 & & & & & & & & & & & & & & & \\
1 & 2 & 3 & 4 & 5 & 6 & 7 & 8 & 9 & 10 & 11 & 12 & 13 & 14 & 15 & 16\\
 & & & & & & & & & & & & & & & \\
\hline
2 & 1 & \cellcolor{mygray}{4} & \cellcolor{mygray}{3} & \cellcolor{mygray}{6} & \cellcolor{mygray}{5} & 8 & 7 & 10 & 9 & \cellcolor{mygray}{12} & \cellcolor{mygray}{11} & \cellcolor{mygray}{14} & \cellcolor{mygray}{13} & 16 & 15\\
\hline
3 & 4 & 1 & 2 & \cellcolor{mygray}{7} & \cellcolor{mygray}{8} & \cellcolor{mygray}{5} & \cellcolor{mygray}{6} &
11 & 12 & 9 & 10 & \cellcolor{mygray}{15} & \cellcolor{mygray}{16} & \cellcolor{mygray}{13} & \cellcolor{mygray}{14}\\
\hline
5 & 6 & 7 & 8 & 1 & 2 & 3 & 4 & 13 & 14 & 15 & 16 & 9 & 10 & 11 & 12\\
\hline
8 & \cellcolor{mygray}{7} & 6 & \cellcolor{mygray}{5} & \cellcolor{mygray}{4} & 3 & \cellcolor{mygray}{2} & 1 &
16 & \cellcolor{mygray}{15} & 14 & \cellcolor{mygray}{13} & \cellcolor{mygray}{12} & 11 & \cellcolor{mygray}{10} & 9\\
\hline
9 & \cellcolor{mygray}{10} & \cellcolor{mygray}{11} & 12 & \cellcolor{mygray}{13} & 14 & 15 & \cellcolor{mygray}{16} &
1 & \cellcolor{mygray}{2} & \cellcolor{mygray}{3} & 4 & \cellcolor{mygray}{5} & 6 & 7 & \cellcolor{mygray}{8}\\
\hline
12 & 11 & \cellcolor{mygray}{10} & \cellcolor{mygray}{9} & \cellcolor{mygray}{16} & \cellcolor{mygray}{15} & 14 & 13 & \cellcolor{mygray}{4} & \cellcolor{mygray}{3} & 2 & 1 & 8 & 7 & \cellcolor{mygray}{6} & \cellcolor{mygray}{5}\\
\hline
14 & 13 & 16 & 15 & \cellcolor{mygray}{10} & \cellcolor{mygray}{9} & \cellcolor{mygray}{12} & \cellcolor{mygray}{11} &
\cellcolor{mygray}{6} & \cellcolor{mygray}{5} & \cellcolor{mygray}{8} & \cellcolor{mygray}{7} & 2 & 1 & 4 & 3\\
\hline
15 & \cellcolor{mygray}{16} & 13 & \cellcolor{mygray}{14} & \cellcolor{mygray}{11} & 12 & \cellcolor{mygray}{9} & 10 &
\cellcolor{mygray}{7} & 8 & \cellcolor{mygray}{5} & 6 & 3 & \cellcolor{mygray}{4} & 1 & \cellcolor{mygray}{2}\\
\hline
\end{tabular}\\
\centerline{the table of $[r+8,16s,16n]$}
\end{center}

As the table of $[r+4,8s,8n]$ and the table of $[r+2,4s,4n]$ are both subtables of that of $[r+8,16s,16n]$, we just explain the last table:
\begin{enumerate}\itemsep=0pt
\item Every cell in the f\/irst row of a table which stands for the f\/irst row of the corresponding admissible matrix is the copy of the f\/irst row of~$A$.
\item Every cell in the second row of a table which stands for the rows from second to $r$-th of the corresponding admissible matrix is the copy of the rows from second to $r$-th of~$A$.
\item Every cell in other rows of a table which stand for the $k$-th rows ($k \geq r+1$) of the corresponding admissible matrix is the copy of the f\/irst row of~$A$.
\item For every cell, the number means the levels and the color means the signs.
\end{enumerate}
Using the same discussion of Remark \ref{discussion of judgement}, it is easy to verify that (1)--(3) of~(\ref{eq:admissible matrix}) hold. For the verif\/ication of~(4) of~(\ref{eq:admissible matrix}), it is enough to consider the 8 added rows and verify those entries which are in the same level. Then the admissibility follows by a direct and simple computation.

\subsection*{Acknowledgements}
This research was supported by NSFC 11471186 and NSFC 11571199.

\pdfbookmark[1]{References}{ref}
\LastPageEnding

\end{document}